\documentclass[12pt,reqno]{amsart}

\usepackage[arrow,matrix,curve]{xy}

\usepackage[dvips]{graphicx} 

\usepackage{amssymb, latexsym, amsmath, amscd, array, hyperref
}
\usepackage{breakurl}   %this allows for url breaks in pdf

\usepackage{enumitem} 
\theoremstyle{definition}

\numberwithin{equation}{section}
\newcommand\N {{\mathbb N}} 

\newcommand\R {{\mathbb R}}

\newcommand\sh{\mathbf{sh}}

\newcommand\astr{{{}^{\ast}\hspace{-0.5pt}\R}}

\newcommand\ZF{\textbf{ZF}}

\newcommand\ZFC{\textbf{ZFC}}

\newcommand\AC{\textbf{AC}}

\newcommand\ACC{\textbf{ACC}}

\newcommand\IST{\textbf{IST}}

\newcommand\SIA{\textbf{SIA}}

\newcommand\LV{\textbf{LV}}

\newcommand\FL{\textbf{FL}}

\numberwithin{equation}{section}

\author[E. Bottazzi]{Emanuele Bottazzi} \address{E. Bottazzi,
Department of Civil Engineering and Architecture of the University of
Pavia, Via Adolfo Ferrata, 3, 27100 Pavia, Italy}
\email{emanuele.bottazzi@unipv.it}

\author[V. Kanovei]{Vladimir Kanovei} \address{V. Kanovei, IPPI RAS,
Moscow, Russia}\email{kanovei@googlemail.com}

\author[M. Katz]{Mikhail G. Katz}\address{M. Katz, Department of
Mathematics, Bar Ilan University, Ramat Gan 52900
Israel}\email{katzmik@macs.biu.ac.il}

\author[T. Mormann]{Thomas Mormann}\address{T. Mormann, Department of
Logic and Philosophy of Science, Univ.  of the Basque Country UPV/EHU,
20080 Donostia San Sebastian, Spain}
\email{thomasarnold.mormann@ehu.eus}

\author[D. Sherry]{David Sherry} \address{D. Sherry, Department of
Philosophy, Northern Arizona University, Flagstaff, AZ 86011, US}
\email{David.Sherry@nau.edu}

\begin{document}

%\doublespacing
\thispagestyle{empty}

%\huge

\title[On mathematical realism and applicability of hyperreals] {On
mathematical realism and the applicability of hyperreals}

\begin{abstract}
We argue that Robinson's hyperreals have just as much claim to
applicability as the garden variety reals.

In a recent text, Easwaran and Towsner (ET) analyze the applicability
of mathematical techniques in the sciences, and introduce a
distinction between techniques that are \emph{applicable} and those
that are merely \emph{instrumental}.  Unfortunately the authors have
not shown that their distinction is a clear and fruitful one, as the
examples they provide are superficial and unconvincing.  Moreover,
their analysis is vitiated by a reliance on a naive version of object
realism which has long been abandoned by most philosophical realists
in favor of truth-value realism.  ET's argument against the
applicability of hyperreals based on automorphisms of hyperreal models
involves massaging the evidence and is similarly unconvincing.  The
purpose of the ET text is to argue that Robinson's infinitesimal
analysis is merely instrumental rather than applicable.  Yet in spite
of Robinson's techniques being applied in physics, probability, and
economics (see e.g., \cite[Chapter\;IX]{Ro66}, \cite{Al86},
\cite{Va07}, \cite{Lo15}), ET don't bother to provide a meaningful
analysis of even a single case in which these techniques are used.
Instead, ET produce page after page of speculations mainly imitating
Connesian chimera-type arguments `from first principles' against
Robinson.  In an earlier paper Easwaran endorsed real applicability of
the~$\sigma$-additivity of measures, whereas the ET text rejects real
applicability of the axiom of choice, voicing a preference for \ZF.
Since it is consistent with \ZF\;that the Lebesgue measure is
\emph{not}~$\sigma$-additive, Easwaran is thereby walking back his
earlier endorsement.  We note a related inaccuracy in the textbook
\emph{Measure Theory} by Paul Halmos.  ET's arguments are unacceptable
to mathematicians because they ignore a large body of applications of
infinitesimals in science, and massage the evidence of some crucial
mathematical details to conform with their philosophical conclusions.

\end{abstract}

\maketitle

\tableofcontents

\section{Varieties of realism}
\label{s1b}

In a text entitled ``Realism in Mathematics: The Case of the
Hyperreals'' \cite{ET}, K. Easwaran and H. Towsner (henceforth ET)
deal with the following two general questions related to the issue of
mathematical realism:
\begin{enumerate}
\item
Which mathematical claims can be taken as meaningful and true within
mathematics?
\item
Which mathematical ideas can be taken as applying to the physical
world as part of scientific theory?
\end{enumerate}
For the special case of hyperreals~$\astr$ in the sense of Robinson's
framework for analysis with infinitesimals and related approaches,
these questions are answered by ET as follows:
\begin{enumerate}
[label={(\theenumi)*}]
\item
\label{1*}
The hyperreals~$\astr$ are meaningful mathematically since their
existence can be ensured within mathematics by the axiom of
choice\;(\AC).  ET consider \AC\;to be true.
\item
\label{2*}
In contrast to the reals~$\R$, the hyperreals~$\astr$ have no real
application in science as part of a scientific theory. Hyperreals may
only be used instrumentally, as `computational tools,' but they do not
correspond to anything that exists in the real world.  Real numbers
correspond directly to objects existing in the real world and
therefore are to be considered as applicable.
\end{enumerate}

ET consider item~\ref{1*} as a claim that supports mathematical
realism since it ensures a certain kind of existence of certain
mathematical entities.  Namely, hyperreals exist within the universe
of mathematical entities.  Complementarily, ET conceive item~\ref{2*}
as an anti-realist claim, since it denies a certain kind of existence
to~$\astr$.  More precisely, hyperreals do not exist in the real
world, since the hyperreals are not part of a (true) scientific theory
of the real world.  Only real numbers exist, since they are part of
(true) scientific theories of the real world.

According to ET, this `mixed' character of their account of
(mathematical) realism shows that ET's brand of realism is more
balanced than other purely `realist' or purely `anti-realist'
accounts.

In this paper we argue that items~\ref{1*} and \ref{2*} actually
provide evidence for something quite different, namely, that ET
subscribe to an outdated and simplistic account of mathematical
realism (and of scientific realism as well) that fails to meet the
criteria of current literature in contemporary philosophy of
mathematics and philosophy of science.  We will first argue for these
claims for mathematical realism.  Then we show that ET not only miss
the point for mathematical realism but for scientific realism as well.

More precisely, ET subscribe to an \emph{ontological realism} in
mathematics that is exclusively concerned with the question of the
\emph{existence} of mathematical objects, be it the existence within
mathematics or outside mathematics.  For some time, this ontological
question (be it in mathematics or in science) is no longer the only
game in town.

\subsection{Dummett, Shapiro,  Putnam}
\label{s11}

With respect to mathematical realism it is generally admitted that
there are at least two different realist themes in contemporary
philosophy of mathematics. According to the first, mathematical
objects exist independently of mathematicians, their minds, languages,
and so on.  This is the only version of realism ET deal with. It may
be succinctly characterized as an ontological or \emph{object
realism}.  A unique feature of ET's approach is their quest for
``numbers that correspond to the world,'' not found in, and a likely
embarrassment for, any object realist such as Colyvan \cite{Co98}.

The second theme in contemporary mathematical realism centers around
the thesis that mathematical \emph{statements} have objective
truth-values independent of the minds, languages, conventions, and so
forth, of mathematicians.  This realism may be called truth-value
realism or \emph{objectivity realism}.  It can be characterized as a
realism that centers around the objectivity of mathematical discourse.
According to it, the interesting and important questions of philosophy
of mathematics are not over mathematical objects, but over the
objectivity of mathematics (Shapiro \cite{Sh97}, 1997, p.\;37).  

Truth-value realism does not center around an \emph{ontological} view
thesis.  Although truth-value realism claims that many mathematical
statements have unique and objective truth-values, it is not committed
to a distinctively platonist idea that such truth-values are to be
explained in terms of an ontology of mathematical objects.  A related
challenge to the uniqueness of ontology of mathematical objects was
developed by Benacerraf \cite{Be65}.

Philosophers of quite different orientations have argued for
objectivity realism in mathematics. Even if they disagree widely on
the various issues related to the existence of mathematical objects,
virtually all converge in the assertion that mathematics is the
\emph{objective} science par excellence.  Therefore they consider it a
central task of philosophy of mathematics to elucidate what exactly is
meant by this assertion.  As Michael Dummett put it succinctly: ``What
is important is not the existence of mathematical objects but the
objectivity of mathematical statements'' \cite[p.\;508]{Du81}.

A prominent example of a philosopher of mathematics who argued
vigorously for a truth-value realism in mathematics is Hilary Putnam.
Although Putnam was notorious for changing his philosophical views a
number of times, with respect to the issue of a non-ontological
realism in mathematics his convictions have remained remarkably
stable.  For forty years in the course of his entire philosophical
career he insisted on a non-ontological realism.  Already in his
\emph{What is mathematical Truth?} he put forward the thesis that
``The question of realism is the question of the objectivity of
mathematics and not the question of the existence of mathematical
objects'' (Putnam\;\cite{Pu75}, 1975, p.\;70). 

In ``Indispensability Arguments in the Philosophy of Mathematics''
(Putnam \cite{Pu12}, 2012) he insisted that his well-known
\emph{indispensability argument} should be understood as an argument
for the objectivity of mathematics, and not as an argument for a
platonist interpretation of mathematics.  The point is reiterated in
his posthumously published essay in 2016: 
\begin{quote}
[O]ne does not have to `buy' Platonist epistemology to be a realist in
the philosophy of mathematics. The modal logical picture shows that
one doesn't have to `buy' Platonist ontology either.''
(Putnam 2016 \cite{Pu16}, p.\;345)
\end{quote}

\subsection{Maddy on varieties of realism}
\label{mad}

Some years ago, Penelope Maddy in her book \emph{Defending the Axioms.
On the Philosophical Foundations of Set Theory} (Maddy \cite{Ma11},
2011) took up the issue of mathematical realism.  In \emph{Defending
the Axioms} we find a detailed and balanced discussion of the two
themes of \emph{object realism} versus \emph{objectivity realism} in
mathematics as presented by Shapiro, Dummett, Putnam and others;
cf.\;Maddy (2011, chapter\;(V.1)).

ET refer to Maddy as well but in a rather unfortunate way. They do not
mention Maddy (2011), but seek to use some of her earlier work to
render plausible their distinction between real application and merely
instrumental use.  ET write: 
\begin{quote}
As Penelope Maddy says about physics, `its pages are littered with
applications of mathematics that are expressly understood not to be
literally true: e.g., the analysis of water waves by assuming the
water to be infinitely deep or the treatment of matter as continuous
in fluid dynamics or the representation of energy as a continuously
varying quantity.' [Maddy, 1992, p. 281]'' (Maddy as quoted in
\cite[p.\;5]{ET})
\end{quote}
Yet, as (Maddy 2011) shows, she is no longer interested in arguing for
the existence of mathematical objects analogous to the existence of
``planets and atoms and giraffes, independently of the human mind''
\cite[p.\;1]{ET}, as ET put it.  Quite generally, the issue of which
numbers `correspond to the world' is not one that interests many
philosophers of mathematical realism today.

In her book Maddy is mainly interested in the objectivity of set
theory.  Even if the question of set-theoretic axioms and their
ontological status were to be settled to the satisfaction of
philosophers, and even if all objects of mathematics could be
reconstructed as sets, there is no reason to assume that thereby all
problems of philosophy of mathematics would be solved.  One may ask:
``Do groups, Riemannian manifolds, Hilbert spaces objectively exist?''
For instance, did groups as mathematical objects come into being in
the second half of the 19th century, or did they exist since times
immemorial as may be suggested by the example of finite symmetry
groups of geometrical figures (see e.g., Wussing \cite{Wu69})?

It seems plausible that philosophy of mathematics is more than
philosophy of set theory.  The objectivity of mathematics should not
be discussed in the overly abstract realm of set theory alone.
Problems of mathematical realism do not only arise for sets, even if,
arguably, every mathematical object could be reconstructed as a
structured set.  Indeed, dealing with the question `Are there groups?'
or `Are there Riemannian manifolds?' instead of the analogous, rather
wornout questions `Are there sets?' or `Are there numbers?' has some
advantages, and not only didactic ones.

\subsection{On the dialectics of tool and object}

If science is a \emph{Werdefaktum} (`Fact in Becoming' as the
neo-Kantians such as Cassirer claimed) the objectivity of mathematical
discourse may have to be conceived as a `fact in becoming' as well,
depending on the historically evolving practice of mathematical
discourse (see \cite[Section\;2.1]{13h}).  The mathematical practice
may change the nature of mathematical objects: they may change from
mere tools to full-fledged objects.  From this perspective also the
history of mathematics (and science) becomes important for the issue
of mathematical realism.

\subsection{Maddy on history of mathematics and science}
\label{s14}

(\emph{Defending the Axioms}, pp.\;27--29).  One clear moral for
mathematics in application is that we are not in fact uncovering the
underlying mathematical structures realized in the world.  Rather, we
are constructing abstract mathematical models and endeavoring to make
true assertions about the ways in which they do or do not correspond
to physical facts.  There are rare cases where such a correspondence
is something like isomorphism, as for elementary arithmetic and simple
combinatorics and there are probably others, like the use of finite
group theory to describe simple symmetries.  However, in most cases
the correspondence is something more complex.

For the relation between the physical and the mathematical, an
isomorphism is the great exception. Moreover, ET ignore the most
essential aspects of mathematical conceptualisation.  The
correspondence between the empirical and the mathematical (ideal) is
more complex.  Cassirer drew the conclusion that not the existence of
objects, but objectivity of the method is important; see further in
Section~\ref{s15b}.
 
ET rely on an outdated model of the relation between mathematics and
physics proposed already by Newton and Galileo, namely, that
mathematics can and does discover the mathematical structure of the
world, or better, that the world has a unique mathematical structure.
Such a view is related to the well-known Cantor--Dedekind postulate
according to which the real mathematical line actual describes the
physical line, a position untenable in view of contemporary knowledge
in physics (see Sections~\ref{s1} and \ref{s34}).

Since the application of mathematical concepts typically amounts to an
idealisation, mathematical concepts have no direct correspondence in
the `real world'.  This holds even for a seemingly elementary concept
such as the natural numbers~$\N$.  After all, it is an essential
feature of~$\N$ that it is an infinite structure, and clearly, even
the `small' infinity of~$\N$ has no direct correspondence in the `real
world'; see further in Section~\ref{s27}.

\subsection{Cassirer on object and objectivity}
\label{s15b}

Maddy and most other contemporary philosophers of mathematics, having
been educated in the Anglo-Saxon tradition of analytic philosophy,
tend to attribute the merit of having shifted the attention of
philosophers of mathematics from questions of object realism to
questions dealing with the objectivity of mathematics to a rather
obscure and unelaborated remark of Georg Kreisel's; see e.g.,
\cite[p.\;115, note\;4]{Ma11}.  However influential Kreisel's remark
may have been in recent analytical philosophy of mathematics, Kreisel
was by no means the first to note the relevance of distinguishing
between a realism concerning mathematical objects and a realism
concerning the objectivity of mathematical discourse.

Already in the first decade of the 20th century, the neo-Kantian
philosopher Ernst Cassirer had elaborated in his \emph{Substance and
Function} (\cite{Ca10}, 1910) a neo-Kantian philosophy of science and
mathematics that emphasized (for both mathematics and natural science)
the importance of distinguishing between two interpretations of
realism, namely, object realism and truth-value (objectivity) realism.

It is not very surprising that Cassirer's philosophy of mathematics
has been virtually ignored in the quarters of analytic philosophy. For
a long time Cassirer has been classified as a partisan of some stripe
of continental idealism, giving analytic philosophers an excuse not to
take his work seriously.

The underlying reason for such a lamentable state of affairs may have
been that in traditional 20th century Anglo-Saxon philosophy there is
a conviction that idealist philosophy on the one hand and serious
science and philosophy of science on the other do not go well
together.  Often, idealism plays the role of a strawman to whom all
the vices are attributed that one wishes to criticize.  In the 21st
century one still finds philosophers such as Susan Haack who propagate
virtually the same caricature of idealism that Scheffler put forward
almost 50 years ago: 
\begin{quote}
An idealist holds that everything there is, is mental: that the world
is a construction out of our, or, in the case of the solipsist, his
own, ideas -- subjective idealism; \ldots\;or that the world is itself
of a mental or spiritual character -- objective idealism, as in
Hegel.  (Haack \cite{Ha02}, 2002, p.\;70)
\end{quote}
Evidently, for Haack idealism is not an option to be taken seriously.
For her, idealism is the bogey of realists.  For some time, however,
such an attitude has turned out to be increasingly untenable due to
the fact that Cassirer's neo-Kantian philosophy of science is
re-evaluated as a serious competitor to the classical analytical
philosophy of science, namely, the logical empiricism of the Vienna
Circle and related groups such as Reichenbach's Berlin group; see
e.g., \cite{Co91}, \cite{Fr99}, \cite{Fr01}.

The starting point of Cassirer's account of \emph{critical idealism}
is the insight that the `object' of scientific knowledge - be it
mathematical or empirical knowledge - should not be conceived as a
kind of Kantian `thing-in-itself' beyond all possible experience.
Rather, 
\begin{quote}
If we determine the object not as an absolute substance beyond all
knowledge, but as the object shaped in progressing experience, we find
that there is no `epistemological gap' to be laboriously spanned by
some authoritative decree of thought, by a `transsubjective command.'
For this object may be called `transcendent' from the standpoint of a
psychological individual; from the standpoint of logic and its supreme
principles, nevertheless it is to be characterized as purely
†´immanent.†¡ It remains strictly within the sphere, which these
principles determine and limit, especially the universal principles of
mathematical and scientific knowledge. This simple thought alone
constitutes the kernel of critical `idealism'.''  (\emph{Substance and
Function} \cite{Ca10}, 1910, p.\;297)
\end{quote}
For Cassirer, knowledge never starts with well-determined objects that
are `given' to a cognizing subject.  Rather, objects are the result of
an objectifying process:
\begin{quote}
To know a content means to make it an object by raising it out of the
mere status of givenness and granting it a certain logical constancy
and necessity. Thus we do not know `objects' as if they were already
independently determined and given as objects, - but we know
objectively, by producing certain limitations and by fixating certain
permanent elements and connections within the uniform flow of
experience.  The concept of the object in this sense constitutes no
ultimate limit of knowledge, but is rather the fundamental instrument,
by which all that has become its permanent possession is expressed and
established.  The object marks the logical possession of knowledge,
and not a dark beyond forever removed from knowledge.''
(\emph{Substance and Function}, p.\;303)
\end{quote}
Perhaps one may say that for Cassirer an `object' is generated by a
scientific method of objectification.  Thus (valid) scientific
concepts do not aim to produce `copies' of pre-existing objects;
rather, 
\begin{quote} Scientific \ldots\;concepts are valid, not in
that they copy a fixed, given being, but in so far as they contain a
plan for possible constructions of unity, which must be progressively
verified in practice, in application to the empirical material.  But
the instrument, that leads to the unity and thus to the truth of
thought, must be in itself fixed and secure. \ldots{} \emph{We need,
not the objectivity of absolute things, but rather the objective
determinateness of the method of experience}.''  (\emph{Substance and
Function}, p.\;322; emphasis added)
\end{quote}

This summary of Cassirer's position may suffice to provide evidence
that his `critical idealism' has little in common with the simplistic
caricature of idealism as found in Haack's description of this
philosophical current in \cite[p.\;70]{Ha02}.  Moreover, Cassirer
offers us a more elaborate account of the interplay of objects and
objectivity in the practice of science than Kreisel's succinct remark
that bluntly asserts that for mathematics its objectivity is more
important than the existence of its objects.  Cassirer's critical
idealism treats mathematical and scientific realism on a par: for both
areas he argues for a truth-value realism as opposed to an ontological
object-centered realism. This should be considered as a virtue since
it enables us to overcome the boundaries between mathematics and the
natural sciences that in the age of an ever more mathematicized
science are becoming ever more artificial and obsolete.

\section{Applicable vs instrumental}
\label{s1}

As we already mentioned in Section~\ref{s1b}, ET pursue a distinction
between mathematical ideas meaningful within mathematics and those
meaningful within the sciences.  Section~\ref{s1b} dealt with
philosophical shortcomings of their approach.  In this section we will
focus on the mathematical shortcomings of their analysis.  ET focus on
\begin{quote}
the question of which mathematical claims can be taken as meaningful
and true \emph{within mathematics}, and which mathematical ideas can
be taken as \emph{applying to the physical world} as part of a
scientific theory.''  \cite[p.\;1]{ET} (emphasis added)
\end{quote}
They also pursue a related distinction between mathematical techniques
that can be described as \emph{applicable} (such as Newtonian
mechanics or quantum mechanics) and ones that are merely
\emph{instrumentally} useful in modeling.  ET take the superior
\emph{applicable} techniques to ``accurately describe'' (see
Section~\ref{s13}) or ``correctly model'' the world.

They seek to apply such an \emph{applicable} vs \emph{instrumental}
dichotomy to an appraisal of Robinson's framework \cite{Ro66} for
analysis with infinitesimals.  The main thrust of their text is the
claim that real numbers are applicable whereas hyperreal numbers are
merely instrumental.  We argue that the ET argument fails at several
levels, as follows.

\begin{enumerate}
\item
(\emph{Realism \`a l'ancienne}) ET claim that they adopt the
philosophical position of scientific \emph{realism} and pursue the
issue of which numbers ``correspond to the world''(see e.g.,
Section~\ref{s23}).  However, the issue as stated is not a relevant
issue in contemporary trends in realism in the philosophy of
mathematics, as analyzed in Section~\ref{s11}.  The current literature
addresses not the issue of which mathematical \emph{objects}
correspond to the world, but rather the issue of the
\emph{objectivity} of mathematical discourse; see Section~\ref{s21}.
\item (\emph{Representing the world and Cantor--Dedekind}) ET seek to
drive a wedge between the applicability of the reals and that of
hyperreals on account, in their words, of providing ``representation
of part of the world'' (\cite[p.\;9]{ET}), but their stance is
conditioned on the adoption of the Cantor--Dedekind postulate.  The
latter involves an identification of a Cantor--Dedekind real line as
they understood it with a physical line.  Such a stance is untenable
in view of contemporary knowledge in physics; see Sections~\ref{s6}
and \ref{s34}.
\item
(\emph{Elegance \"uber alles}) Well into their argument, ET pull out
of a hat a criterion of the \emph{elegance} of a mathematical
framework (a fine criterion, to be sure; see Section~\ref{s15}).  They
use it to justify a reliance on a standard uncountable \emph{number
system} of cardinality~$\mathfrak{c}$.  But ET fail to apply their
\emph{elegance} criterion consistently when it comes to adopting a
standard \emph{structure} of cardinality~$\mathfrak c$ over the
traditional reals; see Section~\ref{s38}.  ET thus display a double
standard indicating a partiality of their analysis.
\item
(\emph{Physics to the rescue}) ET claim that ``an applicable use of
the hyperreals would require substantial new developments in
physics.''  However, such a claim ignores the existence of numerous
such applications as developed e.g., in Albeverio et al.\;\cite{Al86};
see Section~\ref{s51} for additional sources.
\item
(\emph{Applicability of {\rm \AC}-dependent entities}) ET claim that
their objections to using hyperreals in mathematical models in physics
apply also to other entities whose existence cannot be proved in
Zermelo--Fraenkel set theory (\ZF) without the use of the Axiom of
Choice (\AC).  However, they focus their critique on Robinson's
framework, and ignore the existence of many applications relying on
other entities commonly used in traditional mathematics, whose
existence cannot be proved without exploiting \AC-related principles.
\item
(\emph{Walking back earlier endorsement}) ET's position is at odds
with an earlier endorsement of $\sigma$-additivity by Easwaran.  For a
discussion on the applicability of the Lebesgue measure without
\AC\;and on other consequences of accepting at face value ET's claim
that their objections to the hyperreals apply also to mathematical
entities independent of \ZF, see Sections\;\ref{s44} through
\ref{s510}.
\item
\label{i5} (\emph{Undercutting and rebutting}) We will
exploit a dichotomy of \emph{rebutting} vs \emph{undercutting}
developed by Easwaran in \cite{Ea15}.  Briefly, \emph{rebutting} an
argument involves showing that its conclusions contradict those
reached in other work published in reliable venues, whereas
\emph{undercutting} involves finding gaps in the argument itself.  We
will use the dichotomy to analyze the ET text \cite{ET}.
\end{enumerate}

\subsection{Chimera and dart strategies}
\label{s12b}

To elaborate further on the last item~\eqref{i5}, note that a related
dichotomy contrasts two possible strategies for challenging the
applicability of a particular fragment of mathematics:
\begin{enumerate}
[label={(Co\theenumi)}]
\item
\label{i1}
(the ``chimera'' strategy) Argue from first mathematical principles
that there are built-in shortcomings in the fragment of mathematics
that would disqualify it from applications;
\item
\label{i2}
(the ``dart'' strategy) Analyze a particular attempt to apply the
fragment of mathematics to a specific area of science and show
concretely how it falls short of the mark.
\end{enumerate}

An example of the \emph{chimera} strategy \ref{i1} is Connes' attempt
in \cite[p.\;14]{Co04} to exploit the Solovay model so as to argue
that hyperreals are chimeras.  His attempt is evaluated in
\cite[pp.\;263, 278, 279 and Section\;8.2, p.\;287]{13c}.

Another example of strategy \ref{i1} involves the \emph{basic} form of
Smooth Infinitesimal Analysis (\SIA).  The basic form has only
nilsquare infinitesimals, which may hinder applications of the
\emph{basic} form of \SIA\;in more advanced applications in physics
and geometry where second-order Taylor expansion and second
differentials appear, such as Newtonian mechanics or curvature of
curves in differential geometry.  Note that there are more advanced
versions of \SIA\;possessing nilpotent infinitesimals of higher order
than~$2$ (see~\cite{Ko06}); our purpose here is merely to give an
example illustrating the chimera strategy~\ref{i1}.

An example of the \emph{dart} strategy~\ref{i2} is Connes' attempt in
\cite[p.\;13]{Co04} to find fault with hyperreals in a specific
application involving throwing darts at a target; his attempt is
evaluated in \cite[Section\;8.1, pp.\;286--287]{13c}.

Between 1994 and 2007, Connes used to offer apriori chimera-type
arguments `from first principles' against hyperreals.  In 2013, two
articles appeared analyzing Connesian \emph{chimera} critiques:
Kanovei et al.\;\cite{13c} and Katz--Leichtnam \cite{13d}.

Connes has toned down his approach since then.  His most recent piece
touching on the nature of infinitesimals no longer seeks to contrast
the respective infinitesimals of Connes and Robinson but rather those
of Newton and Leibniz.  Connes claims that to Newton, an infinitesimal
was a variable taking determined values and tending to zero, whereas
to Leibniz, it was a number (we will refrain from commenting on the
historical merits of such a claim).  Thus, Connes writes that his own
\begin{quote}
new set-up immediately provides a natural home for the `infinitesimal
variables': and here the distinction between `variables' and numbers
(in many ways this is where the point of view of Newton is more
efficient than that of Leibniz) is essential.  It is worth quoting
Newton's definition of variables and of infinitesimals, as opposed to
Leibniz: [`]In a certain problem, a variable is the quantity that
takes an infinite number of values which are quite determined by this
problem and are arranged in a definite order['] (Connes \cite{Co18},
2018, p.\;168)
\end{quote}
In sum, the difference between Newton/Connes and Leibniz/Robinson now
boils down to the fact that Connes leaves a null sequence alone (or
more precisely replaces it by a compact operator with the
corresponding spectrum), whereas Robinson carries out a quotient space
construction involving a nontrivial equivalence relation, resulting in
infinitesimals that are numbers and not merely sequences.

Unlike Connes, Easwaran continues his attempts to market Connesian
\emph{chimera} arguments from first principles, aimed against
Robinson.  Since Robinson's infinitesimal analysis requires \AC\; or
some weaker form of it in order to develop the hyperreals, Easwaran
seeks to undermine the legitimacy of Robinson's framework by attacking
that of~\AC\;(see further in Section~\ref{connes}).

This, however, constitutes philosophical opportunism because in an
earlier paper, Easwaran sought on the contrary to defend the real
applicability of countable additivity, which depends on \AC\;just as a
hyperreal field does; see Section~\ref{s44}.

\subsection{ET \emph{vs} Lotka--Volterra}

ET provide an example of a strategy of type~\ref{i2} in the context of
the Lotka--Volterra theory (LV), a model of predator-prey population
dynamics.  They claim that LV is only \emph{instrumental} rather than
applicable, on the grounds that the real number resulting from a
solution of the LV differential equations is in general noninteger,
and a noninteger obviously cannot literally represent the size of a
mammal population.

This is an example of an ET attempt to exhibit a shortcoming of the
fragment of mathematics represented by LV.  We argue that it is not a
successful example in any meaningful sense, since the trivial
modification~$\lfloor\LV\rfloor$ defeats the non-representability
objection; see Section~\ref{s3}.

In fact ET focus almost exclusively on strategy \ref{i1} when it comes
to challenging Robinson's framework.  Namely, ET seek to argue
inapplicability from first mathematical principles (rather than
examining specific applications).  In this sense, they are committing
the same fallacy as Connes before them, who sought to base a
refutation on Solovay models, undefinability, etc.

We argue that strategy \ref{i1} is rarely a legitimate way of refuting
the applicability of a mathematical theory (unless one is dealing with
pseudoscience like Sergeyev's and Rizza's; see \cite{17g},
\cite{Ka18}, \cite{Kau18}, \cite{Sh19}).  If so, the ET effort fails
on two counts: (a)\;it misses the mark on what it takes to challenge a
fragment of mathematics in a meaningful fashion; (b)\;the critique is
grounded in an obsolete version of the philosophy of object realism.

\subsection{Axiom of choice}

In connection with the status of the axiom of choice (\AC), ET write:
\begin{quote}
[W]hether or not all well-formed mathematical statements have a
definite fact, there is a fact about the Axiom of Choice, and the fact
is that \emph{it is true}, as well as all of its consequences.
\cite[p.\;2]{ET} (emphasis added)
\end{quote}
It may be helpful to clarify the context of this remark.  There is a
spectrum of opinions among mathematicians as to the exact status of
\AC.  One of possible positions is to postulate \AC\;as being
\emph{true}, as ET have done.

Easwaran argued in \cite{Ea14} against the use of hyperreal fields on
the grounds of their alleged reliance on \AC: 
\begin{quote}
My claim is \ldots\;that no physical facts could make one of these
infinitesimals rather than another be the credences of a particular
agent.  Although the Axiom of Choice guarantees that such
hyperreal-valued functions exist, and although these functions are
quite useful to talk about in mathematical contexts, they have
mathematical structure that goes beyond that of credences.  None of
this rules out a certain \emph{instrumental} use of hyperreals.
\cite[p.\;33]{Ea14} (emphasis added)
\end{quote}
We note that declaring \AC\;to be \emph{true} may be a good
instrumental attitude while doing professional mathematics, but
relying on such an assumption while doing professional philosophy
tends to \emph{undercut} the authors' credibility, because it confuses
syntax and semantics.

While ET adopt the position of mathematical realism, theirs is
presumably intended to be a sophisticated kind of realism publishable
in respectable philosophical venues like \emph{Philosophical Review}
where~\cite{Ea14} was published.  A philosophically sound position
(whether realist or not) would have to account for the difference
between theory and model.  Meanwhile, \AC\;is part of the former and
therefore can't be declared \emph{true} in any but an instrumental
sense.  Philosophical realists like Putnam and Maddy don't subscribe
to naive notions of the kind involved in declaring this or that axiom
to be true, but rather argue for the objectivity of mathematical
discourse; see Sections~\ref{mad} and~\ref{s21}.  Furthermore, what ET
write here about \AC\;is at odds with what Easwaran wrote elsewhere
about the real applicability of~$\sigma$-additivity; see
Section~\ref{s44}.

\subsection{Accurately describing the world} 
\label{s13}

ET write: 
\begin{quote}
The question is about the extent to which the theories we build using
various purported mathematical entities \emph{accurately} describe the
world, or are merely instrumentally useful in modeling the world and
making observational predictions.   \cite[p.\;2]{ET} (emphasis added)
\end{quote}
This passage is found in their Section\;1.2 entitled
\emph{Applicability}.  The unique reference ET provide in this section
is Chakravartty \cite{Ch17}.

However, Chakravartty's Stanford entry is about \emph{scientific
realism}.  For him, scientific realism means realism concerning the
empirical sciences, i.e., physics, chemistry, biology etc.  He
explicitly asserts: ``Scientific realism is a realism about whatever
is described by our best scientific theories.''  In passing, he
mentions other kinds of realisms, for instance `external world
realism,' `sense datum realism,' `mathematical realism.'  None of
these varieties of realisms are treated by Chakravartty. Thus, none of
the terms used by Chakravartty in the Stanford article are directly
applicable to the realism ET are interested in, namely, mathematical
realism.  Thus, the ET reference to Chakravartty's article is
gratuitous.

Furthermore, the ET quest for aspects of mathematical models that
would ``accurately describe the world'' is arguably a quixotic one
(see Section~\ref{s14}).  We will interpret the ET position charitably
as a search for a dichotomy between \emph{applicable} (rather than
\emph{accurate}) fragments and merely \emph{instrumental} fragments.
However, ET don't provide any convincing examples to motivate their
dichotomy.  Their example of mammal-counting (see Section~\ref{s3})
fails to deliver on its promise.  The absence of well-motivated
examples \emph{undercuts} the seriousness of their thesis.

\subsection{$\lfloor \text{Lotka--Volterra} \rfloor$} 
\label{s3}

In connection with the study of predator-prey population dynamics, ET
write: 
\begin{quote}
As a simple example, the structure of the natural numbers seems to
accurately represent counts of predators and prey (at least, when
talking about mammals or birds or other macroscopic animals with clear
individuation). But while the real-valued differential equations of
the Lotka--Volterra model can often be useful in predicting or
understanding the ways the predator and prey populations will change,
the infinite precision of various non-integer values that show up are
not taken to represent the actual numbers of predators or prey in the
ecosystem.  \cite [p.\;2] {ET}
\end{quote}
Here ET seek to establish a contrast between mammal-counting and
Lotka--Volterra modeling of a mammal population in an ecosystem.
However, their pastoral example is unconvincing because one can simply
modify the outcome of the Lotka--Volterra (LV) solution by applying
the integer part (floor) function.  With respect to the resulting
model, that we will denote~$\lfloor\LV\rfloor$, one achieves both of
the following objectives:
\begin{enumerate}
\item
The modified model~$\lfloor \LV\rfloor$ yields an integer for
an answer; and
\item
the mathematical essence of the LV model is unchanged in passing from
LV to\;$\lfloor \LV \rfloor$.
\end{enumerate}
The modification we introduced in no way affects the issue of the
alleged inadequacy of the Lotka--Volterra model in describing ``actual
numbers of predators or prey in the ecosystem" as ET put it.  This
\emph{undercuts} their argument since it highlights ET's failure to
exhibit any meaningful shortcoming on the part of the Lotka--Volterra
model that would make it \emph{instrumental} rather than
\emph{applicable}.  ET return to the LV theme later in the article:
\begin{quote}
[B]iologists often represent populations of predator-prey systems with
the Lotka--Volterra differential equations, using real numbers (rather
than integers) to count populations.  \cite[p.\;5]{ET}
\end{quote}
This recurring example is not helpful if their goal is to convince the
reader that they propose a meaningful dichotomy, because the modified
$\lfloor \LV \rfloor$ solution does have integer values, defeating
their pastoral distinction.

\subsection{Vague predicates}

While it is true that the number of mammals in a given ecosystem or
habitat will be a natural number, the number assigned can not be
warranted as an accurate description.  Indeed, a predicate like ``in a
given habitat'' is not an exact predicate, since it admits borderline
cases.  

The observer doing the counting has to make conventional decisions
about which individuals to count as satisfying the predicate.  His
decisions constitute the first step of idealisation.  Animals move in
and out of habitats just as people move in and out of a crowd.  When
one starts counting species, as ecologists do, things get more
complicated.  Thus, idealisation is involved regardless of
whether~$\N$ or~$\R$ is used as a basic number system.

The question of accurate modeling of population dynamics is addressed
in a more subtle manner than what is found in~\cite{ET} in both
undergraduate textbooks \cite{Lo18} and research articles \cite{Ca12}.

\subsection{Are infinitesimals added to $\R$ or found inside $\R$?}
\label{s4}

In connection with number systems used in scientific modeling, ET
write: ``The hyperreals are an extension of the real numbers to
include infinitesimal numbers, etc.''  \cite [p.\;2] {ET}.

Such a viewpoint is a common one and accurately describes many
applications of infinitesimals.  Significantly, ET's philosophical
objections are contingent upon such an \emph{extension} view.
However, there is another approach that does not view the theory
incorporating very small numbers as such an \emph{extension}, namely
Edward Nelson's approach \cite{Ne77} as summarized in
Section~\ref{s27}.

\subsection{Nelson's view}
\label{s27}

In Nelson's framework, one works within the ordinary real line and
finds numbers that behave like infinitesimals there via a foundational
adjustment.  Such an adjustment involves an enrichment of the language
of set theory through the introduction of a one-place predicate
\textbf{st} called ``standard,'' together with additional axioms
governing its interaction with the axioms of Zermelo--Fraenkel set
theory (\ZFC); see Fletcher et al.\;\cite{17f}, Katz--Kutateladze
\cite{15c}, Lawler~\cite{La11} for details.

The predicate \textbf{st} can be viewed as a formalisation of
Leibniz's distinction between assignable and inassignable numbers; see
\cite{13f}, \cite{17b}.  Arguably, Nelson's formalisation of
mathematical analysis as practiced from Leibniz until Cauchy is more
faithful than the formalisation developed by Weierstrass and his
followers, who had to discard infinitesimals since they were unable to
account for them in a satisfactory fashion.

In Nelson's framework \emph{Internal Set Theory} (\IST), an
infinitesimal~$\epsilon$ is a \emph{real} number such
that~$|\epsilon|$ is smaller than every standard positive real number.
Similarly, a nonstandard natural number~$n\in \N$ is greater than
every standard natural number.

The possibility of viewing infinitesimals as being found within $\R$
itself pulls the rug from under ET's philosophical criticisms anchored
in the \emph{extension} view.  Indeed, whatever argument `from first
principles' ET can put forth against hyperreal fields, it must apply
against the real field, as well, because Nelson's approach takes place
within the real field itself, and is therefore immune to ET's
objections. 

For instance, \emph{since the only order-preserving automorphism of
$\R$ is the identity}, nontrivial automorphisms of $\astr$ in the
\emph{extension} approach are not internal, and therefore don't exist
in Nelson's approach, dissolving ET's objection based on
automorphisms.

In sum, what ET present as arguments from first principles turn out to
depend on technical choices in set-theoretic foundations, namely
whether to work with a language limited to~$\{\in\}$ or the richer
language~$\{\in,\mathbf{st}\}$ of Nelson's \IST, a conservative
extension of \ZFC.

\subsection{Cardinality fallacy}
\label{s5}

ET write in a parenthetical remark:
\begin{quote}
We assume, however, that the language is countable, since we are
concerned with theories which will ultimately be used to describe
physical situations.  \cite[p.\;3]{ET}
\end{quote}
It becomes clear in \cite[Section 3.4, p.\;9]{ET} that when ET mention
countable languages, they refer to countably infinite languages, as
opposed to uncountably infinite ones; see Section~\ref{s38} below for
details.

There is a curious assumption here involving an implied connection
between the countable cardinality on the one hand and what ET refer to
as ``physical situations" on the other.  The bland assumption that
infinity (even countable) can have literal meaning with respect to
``physical situations'' seems unwarranted.  Have ET enumerated
countably many physical entities?  The assumption ties in with the
unique aspect of ET object realism discussed in Section~\ref{s11}.

\subsection{Actual real number}
\label{s6}

ET's view concerning the comparative reality if the reals as compared
to the hyperreals is as follows: 
\begin{quote}
[T]he hyperreals extend the reals with infinite numbers: there is a
real~$\omega$ such that~$r < \omega$ for every actual real number~$r$.
\cite [p.\;3] {ET}
\end{quote}
Note that in this passage, ET refer to~$\omega$ as a ``real number".
Such an approach is consistent with Nelson's viewpoint (see
Section~\ref{s27}), one that ET apparently don't wish to adopt,
\emph{undercutting} their argument.  Furthermore, what exactly is an
``actual real number~$r$"?  The passage seem to suggest that ET
subscribe to the Cantor--Dedekind postulate, involving an
identification of a Cantor--Dedekind real line (as Cantor and Dedekind
understood it) with a physical line; see \cite{17f} for a discussion.
Such a stance is untenable in light of contemporary knowledge in
physics; see Section~\ref{s34}.  Furthermore, ET
\begin{quote}
think there are some features that mathematical models of physical
phenomena need to have to be taken in a realist way, and there are
principled reasons to think that hyperreal models will usually lack
these features.  \cite [pp.\;5--6] {ET}
\end{quote}
However, whatever aspect of idealisation ET identify in hyperreal
fields, will also be present in the real field (see Section~\ref{s27}).
Therefore, while their ``principled reasons" may be valid, they fail
to come to their aid in driving a wedge between applicability of reals
and applicability of hyperreals.

\subsection{Reals are ideal like hyperreals}
\label{s23}

In a similar vein, ET write:
\begin{quote}
We don't deny that hyperreals could figure in models of this sort - we
merely assert that when they are so used, we should recognize them as
idealizations that \emph{don't correspond to the world} in the way
that other parts of the models do.  \cite [p.\;5] {ET} (emphasis
added)
\end{quote}
The ET comment about the idealisation aspect of hyperreal fields is
true enough, but it is equally true about the real field, unless one
adopts the Cantor--Dedekind postulate (see Section~\ref{s6}).  This
observation \emph{undercuts} the ET argument.

\section{Mathematical realism, automorphisms}

The authors' stated pursuit of the philosophy of mathematical realism
fails to meet the criteria of current literature in the field.

\subsection{Maddy to the rescue}
\label{s21}

In a passage quoted in Section~\ref{mad}, ET rely on Maddy's 1992 text
concerning idealisation in mathematical physics.  However, Maddy's
position has evolved significantly since 1992.  Particularly, in 2011
Maddy published \emph{Defending the axioms} that analyzes important
distinctions for the discussion of realism in mathematics, such as
\emph{objectivity} of mathematical discourse versus realism concerning
mathematical \emph{objects}; see \cite[pp.\;115--116]{Ma11}.  The
issue of which numbers ``correspond to the world'' (see
Section~\ref{s23}) is not one that interests many philosophers of
mathematical realism today.

\subsection{The realist/antirealist Potemkin village}

ET write: 
\begin{quote}
An \emph{anti-realist} who claims that all scientific theories are
merely models of this sort, with no clear distinction between the
representational and the fictional parts of the theories, may deny the
\emph{cogency} of the distinction we are interested in.  But if one
accepts a distinction between the realism of counting mammals with
integers and the instrumentalism of counting mammals with real
numbers, then one accepts a distinction of the form we would like to
use.  \cite[p.\;5]{ET} (emphasis added)
\end{quote}
Here ET are attempting to hide behind a figleaf of a traditional
realist vs antirealist dichotomy in an attempt to score a point, but
the figleaf is transparent.  One needn't be an antirealist to doubt
that ET have formulated a ``cogent'' distinction (as they put it)
between applicability and instrumentalism, as we argued in
Section~\ref{s3}.  There may well exist a cogent distinction ``between
the representational and fictional parts of [scientific] theories''
but ET have not given us any cogent reason to ``accept \ldots\;a
distinction between the realism of counting mammals with integers and
the instrumentalism of counting mammals with real numbers.''

\subsection{Automorphisms}
\label{s26}

ET first raise the issue of nontrivial \emph{automorphisms} of~$\astr$
in \cite[p.\;6]{ET}.  They seek to capitalize on the existence of
models of hyperreals with nontrivial automorphisms, in contrast with
the reals that are \emph{rigid}, i.e., lack such automorphisms (see
Section~\ref{s27} for the dissolution of their objection in Nelson's
framework).  The argument is that because of the existence of such
automorphisms, no specific infinitesimal can be taken to have
``physical meaning" any better than its image under an automorphism.

The main fallacy of their argument is that \emph{any} model using an
infinite number, say~$H$, will necessarily involve a certain degree of
arbitrariness, since substituting~$2H$ for~$H$ would typically do just
as well in such a model.%
\footnote{To give an elementary example, in calculus the integral of a
Riemann-integrable function~$f$ on a compact interval~$[a,b]$ can be
computed as follows.  One chooses an infinite hypernatural number~$H$,
partitions~$[a,b]$ into~$H$ subintervals of length~$\Delta
x=\frac{b-a}{H}$, computes the corresponding infinite Riemann sum, and
applies the shadow (standard
part):~$\int_a^bf(x)dx=\sh\left(\sum_{i=1}^{H}f(\xi_i)\Delta{}x\right)$.
Replacing~$H$ by~$2H$ throughout would not affect the value of the
integral.}
No automorphisms are needed to argue this type of indeterminacy.
Their focus on the issue of nontrivial automorphisms misses this basic
point.

Furthermore, if as they declared earlier (see Section\;\ref{s5}), ET
wish to work with countable theories, then the quotient field of
Skolem's non-Archimedean integers (see e.g., \cite[Section\;3.2]{13c})
provides a \emph{rigid} non-Archimedean field as can be proved using
Ehrenfeucht's lemma~\cite{Eh73}.%
\footnote{More precisely, the model~$S$ given by Skolem's
non-Archimedean integers is rigid by Ehrenfeucht's lemma.  Rigidity of
this structure is preserved when passing to the quotient field~$Q(S)$
due to a theorem of Julia Robinson \cite{Ro49} that ensures that~$S$
is definable in~$Q(S)$ by a first-order formula.}
The ET ``automorphisms'' argument is thus based on massaging the
evidence.

As Skolem's integers can be naturally embedded in a hyperreal field
(see \cite[Section~2.4]{14a}), the latter can be viewed as the kind of
elegant ``completion'' of the former that ET allow for the real
completion of countable fields of computable reals (see
Section~\ref{s43}), \emph{undercutting} the ET argument.

\subsection{Physical segments}
\label{s34}

ET write:
\begin{quote}
For the example of the numerical representation of distance, the
physical relation is `longer than' and the operation is the
`concatenation' of two segments \ldots\;Given that `longer than' and
`concatenation' can be applied to any pair of distances, that every
distance can be extended and divided, etc.  \cite [p.\;6]{ET}
\end{quote}
Here ET are postulating both indefinite extension and indefinite
divisibility of physical segments, in line with the Cantor--Dedekind
postulate (see Section~\ref{s6}) but contrary to established physical
theory.  ET assume that they can find ``physical" segments that can be
indefinitely divided.  Such an assumption is dubious and goes contrary
to basic principles of quantum mechanics, where infinite divisibility
breaks down at quantum levels.  Continuum understood as \emph{matter}
is grainy by quantum mechanics, in the sense of violating indefinite
divisibility at scales below Planck.

What about continuum understood as \emph{space-time}?  Einstein field
equation
\[
R_{\mu\nu}-\tfrac{1}{2}R\,g_{\mu\nu}+\Lambda{}g_{\mu\nu}=
\frac{8\pi{}G}{c^4}T_{\mu\nu}
\]
shows via the stress-energy tensor $T_{\mu\nu}$, that the graininess
of matter necessarily affects the nature of the metric and curvature
terms in the left-hand side.  Thus according to general relativity,
the metric and the curvature can have physical meaning only at scales
where matter itself has meaning.  This indicates that the indefinite
divisibility of models of space-time used in relativity theory is
necessarily a mere idealisation as far as very small scales are
concerned.

\subsection{Archimedean circularity}

Concerning the issue of physical representation of real numbers, ET
write:
\begin{quote}
Given that `longer than' and `concatenation' can be applied to any
pair of distances, that every distance can be extended and divided,
that `concatenation' and `longer than' commute appropriately, and the
\emph{Archimedean principle}, one can show that a numerical
representation of the conventional sort chosen must in fact exist,
etc.    \cite[p.\;6]{ET} (emphasis added)
\end{quote}
Here ET seem to postulate that there is an Archimedean principle for
the physical relation ``longer than" but ET don't give a source for
such.  Without a source, what the ET are doing is to \emph{assume} the
conclusion they want to reach, namely that a ``representing'' number
system is necessarily Archimedean because the target is assumed to be
Archimedean, revealing a logical fallacy known as \emph{vicious
circle}.  This observation \emph{undercuts} the ET argument.

ET's fallacy is not uncommon.  In a recent article, Walter Dean
discusses 
\begin{quote}
a tradition within measurement theory which questions whether it is
always appropriate to assume that the mathematical structures employed
as scales for length measurement must be Archimedean.
\cite[p.\;336]{De18}
\end{quote}
In footnote~75, he provides sources in the literature on measurement
theory that have expressed ``similar concerns about the empirical
status of the Archimedean axiom.''

\subsection{Changing the subject}
\label{s29}

ET write: 
\begin{quote}
These examples [of how geometry motivates the introduction of square
roots into the number system] so far only motivate the use of
quadratic extensions of the rationals, but \emph{adequate theorizing}
about physical laws motivates the use of more complete subfields of
the real numbers.  \cite[p.\;8]{ET} (emphasis added)
\end{quote}
In this sentence ET have attempted to change the subject of the
discussion.  Until now ET have spoken about physical reality and what
constitutes ``accurate'' (or \emph{applicable}; see Section~\ref{s1})
representation thereof.  Now they have switched to discussing what
constitutes ``adequate theorizing.''  However, their discussion of
``adequate theorizing" involves an equivocation on the meaning of
\emph{adequate}, since adequate theorizing is not the same as accurate
representation.  Their changed focus becomes transparent by the end of
the same paragraph; see Section~\ref{s15}.

\section{Elegant scientific theorizing}

\subsection{Adequacy and elegance}
\label{s15}

ET write:
\begin{quote}
[A] physical system that used only algebraic numbers could not be
phrased in such an \emph{elegant} way as the differential equations
traditional for Newtonian and later theories.  \cite [p\;8] {ET}
(emphasis added)
\end{quote}
Here ET change the focus of what their qualifier \emph{adequate}
refers to.  The adjective no longer refers to (i) existence of an
accurate fit with physical reality, but rather (ii) the
\emph{elegance} of the scientific theory involved (in this case, the
appropriate ordinary differential equations).  Falling back on the
criterion of elegance undermines their own argument, as we argue in
Section~\ref{s37}.

\subsection{Real world quantities, anyone?}
\label{s43}

ET write:
\begin{quote}
It may be that the full set of real numbers is unnecessary - perhaps
the computable real numbers \ldots\;suffice.  However using a
superstructure causes no harm to the applicability of the theory so
long as the values we attempt to assign to observable quantities
belong to the substructure.   \cite[p.\;8]{ET}
\end{quote}
The logic of their argument forces ET to admit that the real numbers
are merely a convenient idealisation; it fact they acknowledge it
repeatedly (see Section~\ref{s35}).  However, the same is true of the
hyperreals, underscoring the fact that their argument is on shaky
ground.  ET write:
\begin{quote}
[E]ven if some field between, say, the algebraic numbers and the
computable reals is the correct model of length, we can work in the
full model of the reals even though, when assigning values to actual
\emph{real world quantities}, we only ever use values from the
subfield.  \cite[p.\;8]{ET} (emphasis added)
\end{quote}
In this passage, ET admit that some numbers in the model may not have
an ``adequate" referent in their original sense of ``adequate" (see
Section~\ref{s29}).  In other words, some real numbers may not be what
they refer to in this remarkable sentence as ``real world
quantities,'' in a unique form of object realism (see
Section~\ref{s11}).

\subsection{Hypercomputation beyond ``actual quantities''}
\label{s35}

ET write:
\begin{quote}
[E]ven if the computable real numbers (or some other substructure of
the reals) appear to suffice for all current physics, the discovery of
some means of hypercomputation could change the correct choice of
substructure.  \cite[p.\;8]{ET}
\end{quote}
All systems ET mentioned so far as potentially corresponding to ``real
world quantities", such as algebraic and computable numbers, are
countable.  Where do the reals come in then?  ET write: 
\begin{quote}
[W]hen an applicable model is a substructure of a rigid larger
mode[l], that larger model can remain applicable: at worst it includes
theoretical states that can never appear in reality, and which are
therefore never assigned to \emph{actual quantities}.
\cite[p.\;8]{ET} (emphasis added)
\end{quote}
If some real numbers are never assigned the ``actual quantities'' as
ET acknowledge, then it clearly follows that~$\R$ is a convenient
idealisation.

Their comment on rigid models seems to endorse the real field and
imply a criticism of hyperreal fields.  However, hyperreal fields
without automorphisms can also be constructed, making the system
rigid.  This is true both with regard to the usual hyperreal fields
relative to suitable structures of cardinality~$\mathfrak{c}$ (see
Section~\ref{s38}), and with regard to countable models with respect
to countable structures (see Section~\ref{s26}).  ET have again failed
to drive a wedge between the applicability of the reals and the
hyperreals in a meaningful way.

\subsection{More ``actual representation''}

ET write:
\begin{quote}
To justify the claim that the hyperreals are an actual representation
of some part of the world (rather than merely a useful computational
tool), the representation must be unique, or at most involve a small
number of arbitrary choices.  \cite[p.\;9]{ET}
\end{quote}
In the context of the attempt by ET to drive a wedge between reals and
hyperreals, this comment appears to assume that the reals are an
``actual representation of the world'' or that at any rate a suitable
subfield of~$\R$ provides such a representation.  This amounts to
accepting the Cantor--Dedekind postulate; see Section~\ref{s6}.  In
line with our policy of interpreting the ET thesis charitably (see
Section~\ref{s13}), we note that, even if one assumes that what they
seek is an \emph{applicable} theory, hyperreal fields are adequate to
the task; see Section~\ref{s35}.

\subsection{Countable languages}
\label{s37}

ET write:
\begin{quote}
[I]f the continuum hypothesis holds, the hyperreals (viewed as an
ordered field) have many automorphisms which fix the reals.  This
continues to hold in any \emph{countably infinite} expansion of the
language of ordered fields -- for instance, we might want to add many
functions arising as solutions to various differential equations, like
exponential and trigonometric functions, in addition to countably many
units and coordinates.  \cite[p.\;9]{ET}
\end{quote}
Earlier in their text, ET do admit an uncountable number system for
reasons of \emph{elegance}, so as to be able to defend the use of~$\R$
as the basic number system.  But their insistence on trimming the
language to countable size does not deliver the desired
disqualification of non-Archimedean systems, since such countable
systems can be constructed that admit no automorphisms, i.e., are
rigid; see Section~\ref{s26}.

If one insists on keeping all of the uncountably many real numbers,
for similar reasons of \emph{elegance}, symbols for all
functions~$f\colon\N\to\N$ should be included in the structure.  In
such case suitable hyperreal fields do become rigid; see
Section~\ref{s38}.

A standard tool in quantum mechanics is the Hilbert space~$\ell^2$ of
square-summable sequences.  This space has the same cardinality
as~$\mathbb N^{\mathbb N}$, namely~${\mathfrak c}$.  Nonetheless,
quantum mechanics is generally thought of as an \emph{applicable}
theory, \emph{undermining} the ET thesis.

\subsection{Uncountable languages and ridigity}
\label{s38}

ET finally admit the following:
\begin{quote}
In a \emph{uncountable} language, the situation is more complicated --
indeed, in a large enough language, the hyperreals can become rigid
[Enayat, 2006] -- but such an infeasible language is, to say the
least, an unusual setting for a scientific theory.  \cite[p.\;9]{ET}
(emphasis in the original)
\end{quote}
ET's emphasis is on the distinction between countable (see
Section~\ref{s37}) and uncountable languages.  Their emphasis makes it
clear that they view a countably infinite language as a suitable
setting for a scientific theory, whereas an uncountably infinite
language is ``infeasible'' and ``to say the least\;\ldots\;unusual''
to such an end.

Since ET don't reveal the details, it may be helpful to point out that
a standard construction of a hyperreal structure over a countable
index set, together with the language including all~$n$-ary functions
on~$\mathbb N$ and all 1-place predicates on~$\mathbb N$, already
gives a rigid model of~$\astr$.%
\footnote{More precisely, due to the availability of a first-order
definable pairing function in~$(\N,+, \cdot)$ (Cantor's pairing
function), all~$n$-ary functions on~$\N$ can be first-order simulated
once we have all~$1$-place predicates.}
Such a structure has cardinality~$\mathfrak c$.

Such a setting is not an unusual setting for a scientific theory and
on the contrary is the standard setting routinely used in physics; see
e.g., Section~\ref{s37}.

\section{New developments in physics}
\label{s5b}

After making some preliminary remarks concerning the continuum
hypothesis, ET proceed to claim the following:
\begin{quote}
There is a slim road here through which new information could change
our view: new physical discoveries could demonstrate the falseness of
the continuum hypothesis and find some way to uniquely distinguish a
particular, rigid, hyperreal structure which has an observable
physical significance.  Because of this possibility, we cannot claim
that the hyperreals couldn't \emph{possibly} be applicable; we can
only claim an applicable use of the hyperreals would require
substantial new developments in physics.  \cite{ET} (emphasis in the
original)
\end{quote}
ET's claim that ``an applicable use of the hyperreals would require
substantial new developments in physics'' is a non-sequitur even by
their own criterion of \emph{rigidity}, since reasonable rigid
hyperreal systems exist regardless of the status of the continuum
hypothesis; see Section~\ref{s38}.

\subsection{Sources on applications of Robinson's framework}
\label{s51}

The mo\-nograph by Albeverio et al.\;(\cite{Al86}, 1986) contains five
hundred pages of meaningful applications of the hyperreals in physics,
\emph{rebutting} the ET thesis.  ET mention the work~\cite{Al86}
briefly as an ``actual proposal \ldots\;for the use of the hyperreals
in science (as in, for instance, [Albeverio et al., 1986])'' in
\cite{ET19}, but provide no substantive evaluation of its merits.
Many applications to physics, probability theory, stochastic analysis,
mathematical economics, and theoretical ecology appear in the books
\begin{itemize}
\item
Robinson (\cite{Ro66}, 1966, chapter IX),
\item
Nelson (\cite{re}, 1987),
\item
Capi\'nski--Cutland (\cite{Ca95}, 1995),
\item
Arkeryd et al.\;(\cite{Ar97}, 1997), 
\item
Faris (\cite{Fa06}, 2006),
\item
Van den Berg and Neves (\cite{Va07}, 2007),
\item
Loeb and Wolff~(\cite{Lo15}, 2015), 
\item
Lobry (\cite{Lo18}, 2018),
\end{itemize}
and in many other sources.

\subsection{Why elementary equivalence?}

ET acknowledge the possibility ``that some physical quantity has
non-Archimedean behavior'' \cite[p.\;10]{ET} but question the need for
elementary equivalence with~$\R$.  Such a need could be argued as
follows.  If science tells us that some quantity grows exponentially,
we would want ``exponential" to have the usual meaning in the whole
non-Archimedean domain.  In other words, we would want the domain to
be elementarily equivalent to~$\R$ with respect to the properties of
the exponential.

Furthermore, we want elementary equivalence for elegance and there is
no need necessarily to have reasons from physics (even though such
reasons are available).  This is in line with the ET admission of
full~$\R$, rather than using just some countable subfield of it, for
reasons of elegance; see Section~\ref{s15}.

\subsection{Quantifier structure}

ET write:
\begin{quote}
[W]hile the hyperreals are meaningful objects worthy of their own
study, there are other contexts (particularly dealing with real
analysis, and the physical theories using it) where they are mere
tools for avoiding some complexity in dealing with quantifier
structure.  \cite[p.\;12]{ET}
\end{quote}
In mathematical pedagogy particularly at freshman level, one of the
main advantages of the hyperreals is providing a simplification of
such ``quantifier structure'' (see e.g., \cite{17h}).  However, ET's
reference to physical theories makes it clear that they are not
limiting their sweeping ``mere tool'' claim echoing Connes (see
\cite{13c}) to pedagogy.  As such, their claim has little basis.

\subsection{Connes and first mathematical principles}
\label{connes}

ET proceed to relate to the rebuttal as developed in \cite{13c} of
Connes' critique, which involves two separate points:
\begin{enumerate}
\item
the issue of undefinability, where Connes claims to be thoroughly
familiar with the ``Polish school of logic" through a seminar he
participated in, and drops hints related to Solovay models that
indicate that he is talking about a purely mathematical issue that he
claims to be a shortcoming of the hyperreals.
\item
the issue of ``constructiveness" which (unlike Bishop) he interprets
as applicability to physics.  Sanders has analyzed the difference
between Bishop's and Connes' take on ``constructiveness" in
\cite{Sa17}.
\end{enumerate}
ET treat the issue as follows.  They mention Connes' critique and the
rebuttal that appeared in \cite{13c}, and then point out that ``the
real issue'' is not the abstruse one of the details of the
mathematical definitions, but rather an alleged inapplicability in
physics:
\begin{quote}
One major line of debate focuses on the claim that no non-standard
hyperreal can be `named'.  Connes tries to establish this claim by
arguing that a given non-standard integer in a hyperreal field can be
used to generate a non-principal ultrafilter over the natural numbers,
or a non-measurable set of reals.  Since these complex entities seem
to be beyond some limit of complexity for physical beings like us to
grasp, this is said to raise problems for any appeal to hyperreals.
Kanovei, Katz, and Morman dispute Connes' claims about the association
of a non-principal ultrafilter with a given hyperreal field \ldots\;We
claim that \emph{definability is not the real issue} here. These
things are `definable' in the strict mathematical sense, but the
`definitions' don't serve the purpose we need definitions to serve in
\emph{physical models}, of making it possible to uniquely measure
quantities.  \cite[p.\;12]{ET} (emphasis added)
\end{quote}
In this passage, ET are equivocating on the meaning of Connes'
criticism and more precisely conflating two separate issues.  Namely,
as argued in \cite{13c}, Connes places himself on the purely
mathematical plane (the \emph{chimera} track of Section~\ref{s12b})
when he voices his \emph{undefinability} critique.  The article
\cite{13c} argues that his critique is incoherent.  ET misrepresent
the picture by changing the subject to physical applications (the
\emph{dart} track of Section~\ref{s12b}).  Now the \emph{dart}
objection is also an objection Connes formulated, but it is a
different one, and the response to that in \cite{13c} was different,
as well.  Thus the ET criticism of the rebuttal of Connes
in~\cite{13c} involves the fallacy of moving the goalposts.

Connes claims that he can establish from first mathematical
principles, having to do with Solovay's model, that something is amiss
with Robinson's infinitesimals (and therefore shift the focus to
Connes' own noncommutative infinitesimals; see also the analysis given
in~\cite{13d} of Connes' critique).  That is the myth debunked
in~\cite{13c}.  This particular theoretical issue has nothing directly
to do with applications to physics (though Connes seeks to apply his
contention so as to claim that allegedly no such applications are
possible).

The ET approach suffers from the same shortcoming as Connes', in that
they ignore the concrete applications of Robinson's framework as for
instance in~\cite{Al86} (see Section~\ref{s51} for other sources), and
instead attempt to argue from first mathematical principles that
Robinson's framework is not applicable to physics.  Their attempt
misses the target as did the earlier salvos of both
Connes\;\cite{Co04} and Easwaran \cite{Ea14}.

\subsection{Applications,~$\sigma$-additivity of Lebesgue measure}
\label{s44}

The ET text contains the following remarkable passage:
\begin{quote}
Note that our worry in application is about existence proofs, and not
all uses of the Axiom of Choice. [Bascelli et al., 2014] point out
(pp.\;851-2) that the countable additivity of Lebesgue measure
requires the Axiom of Choice to prove. It is consistent with \ZF\;that
the set of all real numbers be a countable union of countable
sets. Rejection of the Axiom of Choice would surely cause problems for
those who want to say that Lebesgue measure is countably additive.
\cite[p.\;13]{ET19}
\end{quote}
Having stated the dilemma, ET attempt to resolve it: 
\begin{quote}
But we don't reject the Axiom of Choice - we reject the realist
applicability of mathematical entities whose existence is independent
of \ZF.  Thus, even if someone were to convince us to accept some
alternate mathematics on which there is a countable collection of
countable sets whose union is all of~$\R$, we would say that in
practice it is safe to assume countable additivity of Lebesgue
measure, because these counterexamples would have no applicability in
practice.  (ibid.)
\end{quote}
Here ET allude to models of set theory where~$\R$ is a countable union
of countable sets, e.g., the Feferman--Levy model (\FL).  The
pertinence of the ET claim that such models ``have no applicability in
practice'' needs to be understood.%
\footnote{The \FL\;model of $\R$ admits a decomposition
\begin{equation}
\label{e51}
\R=\bigcup_{n\in\N} C_n \text{ where each } C_n \text{ is countable}.
\end{equation}
Note that \FL\;cannot support a countably additive Lebesgue measure
(because countable sets are Lebesgue measurable and their measure
is~$0$, while countable additivity and decomposition~\eqref{e51} would
then entail that every subset of~$\R$ is a null set).  One can still
define a Lebesgue measure in \FL, but it will only be finitely
additive (see Section~\ref{s55}).  Thus it is consistent with
\ZF\;that the Lebesgue measure is not~$\sigma$-additive.}
In the passage quoted above, ET attempt to change the subject from
\begin{enumerate}
\item
\label{i1b}
the applicability of~$\sigma$-additivity of Lebesgue measure, to
\item
\label{i2b}
the properties of the Feferman--Levy model.
\end{enumerate}
Their strategy is based on the fallacy of moving the goalposts.  What
is important in applications is not the somewhat paradoxical
decomposition properties of item~\eqref{i2b} (see \cite{17i} for an
interesting consequence), but rather the convenient tool of
$\sigma$-additivity of Lebesgue measure as in item~\eqref{i1b}.  It
may be interesting to note that Lebesgue himself required his measure
to be countably additive; see \cite[p.\;236]{Le02} and
\cite[p.\;122]{Ha01}.  In fact, the endorsement of~$\sigma$-additivity
in (Easwaran \cite{Ea13}, 2013) creates an awkward situation where the
author wants to eat the cake (``reject the realist applicability of''
\AC) and have it, too (namely, assume~$\sigma$-additivity of the
Lebesgue measure).

ET claim to ``point out serious problems for the use of the hyperreals
(and other entities whose existence is proven only using the Axiom of
Choice) in describing the physical world in a real way.''  They also
remark:
\begin{quote}
[O]ther entities dependent on the Axiom of Choice don't lend
themselves so naturally to physical theorizing, but we think the
points we make apply generally.  \cite[page\;1]{ET19}
\end{quote}
Meanwhile, the $\sigma$-additive Lebesgue measure is another entity
dependent on \AC\;that is widely used in physical theorizing (see
Section~\ref{s45}). Thus the ET denial of applicability to hyperreals
would apply equally well to the~$\sigma$-additive Lebesgue measure,
contrary to much evidence that points in the opposite direction,
namely to the usefulness of the~$\sigma$-additive Lebesgue measure in
scientific applications, \emph{rebutting} the ET thesis.

In sum, the real issue is not, as ET imply, whether the decompositions
as in~\eqref{i2b} have ``applicability in practice", but rather
whether $\sigma$-additivity mentioned in~\eqref{i1b} has applicability
in practice.

Recall that one of the applications of the Lebesgue measure is the
definition of the spaces~$L^p$ and of the Sobolev spaces~$W^{k,p}$; we
recall that if~$1\leq p \leq \infty$ these are Banach spaces, and
if~$p=2$ the spaces~$L^2$ and the space~$W^{k,2}=H^k$ are Hilbert
spaces. 

Sobolev spaces are widely used for the study of Partial Differential
Equations, that in turn are applied to the mathematical description of
physical phenomena.  The essential point here is that these
applications depend upon properties of the Sobolev spaces that are
independent of \ZF, and typically require \AC\;for their proof; see
Sections~\ref{s55} and~\ref{s45} for more details.

Faced with the ubiquitous reliance on the axiom of countable choice
(\ACC) in the foundations of analysis and topology, many
mathematicians opt to incorporate \ACC\;(or the stronger axiom of
countable dependent choice) as part of the basic foundational package,
and adopt \ZF+\ACC\;as their philosophical credo.  However, this
option is not available to ET since they ground their \emph{real
applicability} thesis (on behalf of mathematics based on \ZF) in the
allegedly constructive nature of \ZF\;foundations (a claim that would
be even less plausible for \ZF+\ACC\;than for~\ZF).  ET's constructive
track is analyzed in Section~\ref{s45}.

\subsection{Connes \& Katz \emph{versus} ET}

ET seek to oppose Connes and Katz on the issue of scientific
applications:
\begin{quote}
In recent years, there has been much debate about the value of the
hyperreals, with two main views exemplified by Alain Connes and
Mikhail Katz, and various coauthors of each.  \cite[p.\;11]{ET}
\end{quote}
However, the salient point here is that Connes and Katz are on the
same side as \emph{against} ET's denial of applicability of
\AC\;and/or related principles.  Thus, Connes routinely exploits
\AC\;in developing the objects he needs to work with his (Connes')
infinitesimals, such as the Dixmier trace (which is a Connesian
version of integration of noncommuting infinitesimals), and exploits
the \v Cech--Stone compactification~$\beta\N$ of~$\N$ in \cite[ch.~V,
sect.\;6.$\delta$, Def.\;11]{Co94}; for details see \cite{13c}.

\subsection{Measures without \AC}
\label{s55}

Few authors have focused on the study of the properties of Sobolev
spaces in \ZF\;without \AC, i.e., in a setting where the Lebesgue
measure is only finitely additive.

Terry Tao in \cite[Definition\;1.2.2]{Ta11} uses a definition of
Lebesgue measure that works in \ZF.  When extended to \ZFC\;this
definition gives the standard~$\sigma$-additive Lebesgue measure.
Therefore Tao's definition is preferable to definitions that
have~$\sigma$-additivity built into the definition and therefore don't
make sense over \ZF.  In view of the above, the measure as defined by
Tao is not merely a Lebesgue-like measure but arguably the Lebesgue
measure itself.  

We note that the measure theory textbook by Paul Halmos mentions
\AC\;only in the context of the construction of a nonmeasurable set.
Halmos' proof of~$\sigma$-additivity makes no mention of \AC, and is
therefore inaccurate.  The gap is in \cite[p.\;42, line -3]{Ha74},
where \ACC\;is relied upon implicitly.  Further gaps in Halmos (mainly
of a philosophical type) are analyzed in~\cite{16b}.

%The comments above may be useful to dispel an impression that
%apparently lingers in part of the mathematical community that
%``Lebesgue measure as familiarly known does not exist without \AC,''
%as claimed by a referee for the \emph{American Mathematical Monthly}
%in a report on the article \cite{17i} which deals with Lebesgue
%measure in a \ZF\;context.  On the basis of the report, the
%\emph{Monthly} rejected the paper, which was subsequently published in
%\emph{Real Analysis Exchange}.
%

The book \cite{Ra83} contains a thorough study of the properties of
finitely additive measures and of the corresponding~$L^p$ spaces.  The
main drawbacks of using finitely additive measures instead of the
Lebesgue measure are that
\begin{enumerate}
\item
the~$L^p$ spaces are not complete with respect to convergence in
 measure;
\item
their completions, the so-called~$V^p$ spaces, are Banach spaces;
however, it is not mentioned whether~$V^2$ is a Hilbert space.
\end{enumerate}
If~$V^2$ is not a Hilbert space, then many results on~$H^k$ spaces
relying on the inner product given by the Lebesgue integral might not
hold for finitely additive measures.  As a consequence, the study of
relatively simple problems, such as the weak formulation of the
Laplace equation, might not be practical, since the existence of a
weak solution to its Dirichlet problem relies on the Riesz
Representation Theorem for the Sobolev space~$H^1_0$ and on the weak
sequential compactness of this space.  We recall that the Laplace
equation is used as a model of various physical phenomena, since its
solution can be interpreted as the density of a physical quantity in a
state of equilibrium.  Some physical laws, such as Fourier's law of
heat conduction or Ohm's law of electrical conduction, can be
formulated by means of the Laplace equation; see Evans \cite{Ev98} for
further details.

We remark that, by working with the algebra of Borel-coded
subsets~$\mathcal{B}$ of a topological space~$X$ instead of the whole
algebra of Borel subsets, it is possible to refine some finitely
additive measures defined over~$X$ to \emph{Borel-coded measures},
that are~$\sigma$-additive over a suitable
$\sigma$-algebra~$\mathcal{E}\subseteq\mathcal{B}$. The details of the
construction are discussed by Fremlin~\cite[Chapter\;56]{Fr08}; we
refer also to the appendix of Foreman--Wehrung~\cite{Fo91} for a
gentle introduction and further references.  If the underlying
topological space~$X$ is second countable, the spaces~$L^p(X)$, whose
elements are (equivalence classes) of Borel-coded functions whose
$p$-th power is integrable with respect to a codably~$\sigma$-finite
Borel-coded measure, are norm-complete whenever~$1\leq p < \infty$;
under the same hypotheses~$L^2$ is also Hilbert space
\cite[Chapter\;56, pp.\;204, 212]{Fr08}.

Note that this construction applies to the finitely-additive Lebesgue
measure over~$\R$.  Since second countability of~$\R$ can be proved in
\ZF\;alone (see Herrlich~\cite{He06}), even in the Feferman--Levy
model of the real line there is a (Borel-coded) Lebesgue measure which
is~$\sigma$-additive on a nontrivial~$\sigma$-algebra of subsets
of~$\R$ and whose~$L^p$ spaces are complete.  However, to the best of
our knowledge, these spaces of Borel-coded Lebesgue measurable
functions have not yet been successfully applied in the sense
advocated by ET.\, See Section~\ref{s45} for further limitations on
the applications of the~$L^p$ spaces without \AC.

\subsection{ET's worry and the Dirichlet problem}
\label{s45}

In the passage quoted in Section~\ref{s44}, ET express a ``worry''
concerning existence proofs that make use of \AC.  They advocate the
use of ``constructive or computable approximations"
\cite[p.\;13]{ET19} of some results depending on it, such as the
Hahn--Banach theorem.  Their position echoes the more radical
rejection of indirect proofs by constructive mathematicians following
Bishop and Bridges \cite{[1]}, \cite{[2]}.  

We recall that the approach of Bishop and his school is claimed to be
consistent with classical mathematics without \AC, since at its core
it consists in the rejection of the law of excluded middle and of
those principles of classical mathematics that imply it, such as
\AC\;or the Hahn--Banach theorem.  Existence results obtained without
these principles can and have been turned into algorithms; however,
their scope does not yet include some central areas of mathematics
such as the theory of partial differential equations.

Consider for instance the constructive Laplace equation, discussed by
Bridges and McKubre-Jordens in \cite{[3]}.  The authors admit that it
is not always possible to define constructively a solution (since the
Riesz Representation Theorem is not constructively valid for every
function in~$H^1_0$ and since weak sequential compactness does not
entail constructively the existence of a limit).  Thus, they write:
``We prove the (perhaps surprising) result that the existence of
solutions in the general case is an essentially nonconstructive
proposition: there is no algorithm which will actually compute
solutions for arbitrary domains and boundary conditions''
\cite[p.\;1]{[3]}.  In a similar vein, Bridges and McKubre-Jordens
write:
\begin{quote}
In this section we prove that the existence of a weak solution of the
general Dirichlet problem [for the Laplace equation] for a
domain~$\Omega\subset \R^2$ cannot be proved constructively.
\cite[p.\;6]{[3]}
\end{quote}
In fact their article may have been more accurately titled ``Failing
to solve the Dirichlet problem constructively.''

As a consequence, this relatively simple equation is not yet tractable
by means of their approach.  Meanwhile, the non-constructive proof
techniques that are used in classical mathematics to prove the
existence and uniqueness of a solution to the Laplace equation are
applied systematically to other PDEs arising in physics and in
engineering (see Pinchover--Rubinstein \cite{Pi05}).  To the best of
our knowledge, no applicable alternative to the use of \AC\;in the
theory of PDEs has been proposed yet.

\subsection{Applications of finitely additive measures}

Some mathematicians have argued against the exclusive use
of~$\sigma$-additive measures.  For instance, de Finetti suggested
that fair lotteries over infinite sets should be modelled by finitely
additive measures \cite{deFinetti}.  With different motivations,
Nelson proposed an approach to probability theory based solely on
hyperfinitely additive measures in the context of Internal Set Theory
\cite{re}; see also the more recent axiomatic approaches by Benci et
al.\;\cite{nap1}, \cite{nap2}.

The advantages of relinquishing~$\sigma$-additivity are not limited to
probability theory: a well-known result in Robinson's framework
entails that any measure, be it~$\sigma$-additive or simply finitely
additive, can be represented by a hyperfinite counting measure (see
Henson \cite{henson}); in addition, it is possible to require many
degrees of compatibility between the standard measure and its
hyperfinite representative \cite{bbdn}, or even construct a single
hyperfinite counting measure that is simultaneously compatible with
all of the Hausdorff measures \cite{Wa77}.

Hyperfinite counting measures enable one to define functional spaces
that are expressive enough to represent not only the Sobolev spaces
$W^{k,p}$, but also other generalized functions commonly used for the
study of PDEs. For instance, in \cite{bottazzi2019} it is shown that
all linear PDEs and many nonlinear ones can be given an equivalent
nonstandard formulation in the space of \emph{grid functions} of
nonstandard analysis, thus providing a unifying framework for the
study of problems that in standard mathematics require different
approaches. For a more precise statement on the advantages of grid
functions in the theory of PDEs, we refer to \cite{bottazzi2019},
\cite{preprint}, and \cite{thesis}.

\subsection{Effective concepts}
\label{s510}

Brunner et al.\;open their article with the following illuminating
passage:
\begin{quote}
A concept is \emph{effective} in the sense of Sierpi\'nski if it does
not require the axiom of choice \AC.  Here we show by means of
examples that fundamental notions of quantum theory are not effective.
For instance (see Section\;1.2) there is an irreflexive Hilbert
space~$L$, constructed from Russell's socks in the second Fraenkel
model~$\mathcal M_2$.  Hence the very notion of a self-adjoint
operator as an observable of quantum theory may become meaningless
without the axiom of choice.  (Brunner et al.\;\cite[p.\;319]{Br96})
\end{quote}
The authors go on to hedge their bets with the following:
``Nevertheless we identify a nontrivial class of observables, the
intrinsically effective Hamiltonians, which is compatible with~$L$ in
the following sense, etc." (ibid.).  We should note that a typical
physicist is not very interested in limiting the scope of
applicability of mathematical results by introducing foundational
restrictions (such as banning \AC), nor in introducing technical
complications necessitated by such restrictions.  On the contrary, he
seeks constantly to ``push the envelope" by applying mathematical
methods somewhat beyond their `official' domain of applicability.
Well-known examples of such attitudes are Dirac's delta `function' and
the Feynman `integral.'

\section{Conclusion}

Easwaran and Towsner argue that the hyperreal number system of
Robinson's infinitesimal analysis is a good instrumental theory, i.e.,
a theory that is ``useful for proving theorems about the real
numbers,'' but at the same time it is not suitable for the description
of physical phenomena, since the infinitesimals of Robinson's
framework are ``idealisations that don't correspond to the world.''
However, these claims rest upon an outdated conception of mathematical
realism and on the identification of the physical continuum with the
Cantor--Dedekind continuum of real numbers as understood by these
classical authors.  By abandoning the idea that the relation between
mathematics and physical reality must be that of an isomorphism, one
sees that the arguments proposed by ET against the applicability of
the hyperreals, such as the distinction between \emph{instrumental}
and \emph{applicable} theories, are inadequate.

While ET concede that new developments in physics could lead to a
different model of the physical continuum, they seem to ignore the
fact that this approach has already been proposed by many authors. In
fact, there are many applications of Robinson's framework to physics,
economics and other sciences, that are unfortunately not discussed in
any detail by ET, who choose instead to argue that the hyperreals are
inapplicable from first principles and with ad-hoc arguments.

In an attempt to show that their critique is not meant to single out
Robinson's framework, ET suggest that the flaws attributed to the
hyperreal numbers are shared also by other mathematical entities whose
existence cannot be proved in \ZF, being dependent on the axiom of
choice (\AC).  Yet, they quickly conclude that ``other entities
dependent on the Axiom of Choice don't lend themselves so naturally to
physical theorizing'' and ``they won't generally play any role in
application,'' thus focusing their attack on Robinson's
infinitesimals.  However, their stance on the applicability of
mathematical entities dependent upon \AC\;is inconsistent, since E
argued in an earlier publication in favor of~$\sigma$-additivity, a
property independent of \ZF.  Moreover, ET do not address the fact
that many other applicable mathematical theories require some form of
a choice principle.  As an example, in Section~\ref{s5b} we showed
that even simple PDEs become intractable without~$\sigma$-additivity
of the Lebesgue measure.

Thus, if mathematical entities dependent upon \AC\;should cause worry
and ultimately are deemed not truly applicable, then according to ET
many areas of mathematics, such as the theory of PDEs or the formalism
of quantum mechanics, would suffer from the same drawbacks attributed
to the hyperreals. Being that these theories are widely applied also
outside of mathematics, the ET position is hardly defensible.

%In Section~\ref{s5b} we argued that there already exist ``new
%developments in physics'' based on the hyperreals, so that they
%already possess relevant applications to the physical sciences.  We
%also argued that the ET stance on the applicability of mathematical
%entities that depend upon \AC\;is inconsistent, since they reject the
%hyperreals and express ``worry'' concerning the use of mathematical
%entities whose existence depends upon \AC, but E argued in an earlier
%publication in favor of the~$\sigma$-additive Lebesgue measure,
%similarly dependent on \AC.  Moreover, if mathematical entities that
%depend upon \AC\;should cause worry (and ultimately are deemed not
%truly applicable), then according to ET many areas of mathematics,
%such as the theory of PDEs or the formalism of quantum mechanics,
%would suffer from the same drawback.  Yes such a position is hardly
%defensible, since there are many examples of applications of
%mathematical theories dependent upon \AC.  In fact, physicists also
%apply mathematical theories in the making, as in the examples provided
%by Dirac and Feynman (that, once formalized, still rely on \AC).  As
%noted by Karl Svozil, ``Physical and set theoretical entities must be
%operationalized wherever possible.  At the same time, physicists
%should be open to `bizarre' or `mindboggling' new formalisms, which
%need not be operationalizable or testable at the time of their
%creation, but which may successfully lead to novel fields of
%phenomenology and technology'' \cite[p.\;1541]{Sv95}.

\section*{Acknowledgments}

We are grateful to Ali Enayat, Karel Hrbacek, and Claude Lobry for
helpful comments on an earlier version of the article.

\end{document}